\documentstyle[12pt]{article}
 \newcommand{\sect}[1]{\setcounter{equation}{0}\section{#1}}

 \font\tengoth=eufm10 \font\sevengoth=eufm7 \font\fivegoth=eufm5
 \newfam\gothfam
 \textfont\gothfam=\tengoth \scriptfont\gothfam=\sevengoth
   \scriptscriptfont\gothfam=\fivegoth
   \def\goth{\fam\gothfam}    %% Euler Fraktur (math mode only)
 \font\frak=eufm10 scaled\magstep1
 \newcommand{\fra}[1]{\mbox{\frak #1}}

 \def\be{\begin{equation}}
 \def\ee{\end{equation}}
 \def\bea{\begin{eqnarray}}
 \def\eea{\end{eqnarray}}

\def\g{\mbox{${\fra g}$}}

\def\gdps{\mbox{${\fra g}_{\cal P'}^*$ }}
\def\gde{\mbox{${\fra g}_{\cal E}^*$ }}
\def\gdest{\mbox{${\fra g}_{\cal E'(\theta)}^*$ }}

\def\gdr{\mbox{${\fra g}_{\cal R}^{*}$ }}
\def\gddj{\mbox{${\fra g}_{\cal DJ}^{*}$ }}
\def\gddjr{\mbox{${\fra g}_{\cal DJR}^{*}$ }}

 \def\L{\mbox{${\bf L}$}}

\def\faj{\mbox{${\cal F}_{j}$ }}
\def\fae{\mbox{${\cal F}_{\cal E}$ }}

\def\faps{\mbox{${\cal F}_{\cal P'}$ }}

\def\far{\mbox{${\cal F}_{\cal R}$ }}

\def\farol{{\cal F}_{\widetilde{\cal R}(\lambda)}}
\def\fart{\mbox{${\cal F}_{{\cal R}(\theta)}$}}
\def\fapsrol{{\cal F}_{{\cal P'}\widetilde{\cal R}(\lambda)}}
\def\faeab{\mbox{${\cal F}_{\cal E(\alpha,\beta)}$}}
\def\faesab{\mbox{${\cal F}_{\cal E'(\alpha,\beta)}$ }}

\def\mdj{\mbox{$\mu_{\cal DJ}^*$}}
\def\mdjr{\mbox{$\mu_{\cal DJR}^*$}}

\def\mre{\mbox{$\mu_{\cal R}^*$}}
\def\mret{\mbox{$\mu_{\widetilde{\cal R}}^*$}}
\def\mps{\mbox{$\mu_{\cal P'}^*$}}

\def\mes{\mbox{$\mu_{\cal E'}^*$}}

\def\uneosl{\mbox{$U_{\cal E}(sl(3))$}}
\def\uneabsl{\mbox{$U_{\cal E(\alpha , \beta )}(sl(3))$}}
 
\def\unpsosl{\mbox{$U_{\cal P'}(sl(3))$}} 
\def\undjosl{\mbox{$U_{\cal DJ}(sl(3))$}} 
\def\undjrosl{\mbox{$U_{\cal DJR}(sl(3))$ }} 
\def\undjrlsl{\mbox{$U_{\cal DJR(\lambda)}(sl(3))$ }}
\def\undjrtsl{\mbox{$U_{\cal DJR(\theta)}(sl(3))$ }}

\def\dele{\mbox{$\Delta_{\cal E}$}}

\def\delps{\mbox{$\Delta_{\cal P'}$ }}
\def\delpsol{\mbox{$\Delta_{{\cal P'}\widetilde{\cal R}(\lambda)}$}}

 \def\deleabh{\mbox{$\Delta_{\cal E(\alpha , \beta)}\,(H)$}}
 \def\deleaba{\mbox{$\Delta_{\cal E(\alpha , \beta)}\,(A)$}}
 \def\deleabb{\mbox{$\Delta_{\cal E(\alpha , \beta)}\,(B)$}}
 \def\deleabe{\mbox{$\Delta_{\cal E(\alpha , \beta)}\,(E)$}}

 \def\delesabh{\mbox{$\Delta_{\cal E'(\alpha , \beta)}\,(H)$}}
 \def\delesaba{\mbox{$\Delta_{\cal E'(\alpha , \beta)}\,(A)$}}
 \def\delesabb{\mbox{$\Delta_{\cal E'(\alpha , \beta)}\,(B)$}}
 \def\delesabe{\mbox{$\Delta_{\cal E'(\alpha , \beta)}\,(E)$}}

 \def\exs{\mbox{$e^{\sigma}$}}
 \def\exms{\mbox{$e^{-\sigma}$}}
 
 \def\exmds{\mbox{$e^{-2\sigma}$}}

 \def\exdst{\mbox{$e^{2\widetilde{\sigma}}$}}
 \def\exmdst{\mbox{$e^{-2\widetilde{\sigma}}$}}
 
 \def\exmcst{\mbox{$e^{-4\widetilde{\sigma}}$}}
\begin{document}
\begin{center}
 
{\LARGE{\bf{ Extended and Reshetikhin Twists for
$sl(3)$}}}\footnote{This  work has been partially supported by DGES of the 
Ministerio de  Educaci\'on y Cultura of  Espa\~na under Projects
PB95-0719 and  SAB1995-0610, the Junta de Castilla y Le\'on (Espa\~na),
and the  Russian  Foundation for Fundamental Research under grant 
97-01-01152.}\\[3mm]

\vskip1cm
 
{\sc Vladimir D. Lyakhovsky{\footnote {In absence from {Theoretical
Department, Sankt-Petersburg State University; 198904,  St. Petersburg,
Russia.}} and Mariano A. del Olmo} 
\vskip0.5cm
 
{\it  Departamento de  F\'{\i}sica Te\'orica,
 Universidad de Valladolid,  }\\
 {\it E-47011, Valladolid,  Spain}}
\vskip0.25cm

\vskip0.15cm

{email: vladimir@klander.fam.cie.uva.es, 
olmo@fta.uva.es}
 \end{center}
 
\vskip1cm
\centerline{\today}
\vskip1.5cm

\begin{abstract}
The properties of the set ${\cal L}$ of extended jordanian twists 
for algebra $sl(3)$ are studied. Starting from the simplest algebraic 
construction --- the peripheric 
Hopf algebra $U_{{\cal P'}(0,1)}(sl(3))$ ---
we construct explicitly the complete family of extended twisted algebras 
$\{ U_{\cal E(\theta)}(sl(3)) \}$ corresponding to the 
set of 4-dimensional Frobenius subalgebras $\{ {\bf L(\theta)} \}$ 
in $sl(3)$. It is proved that the extended 
twisted algebras with different values of the parameter $\theta$ are 
connected by a special kind of Reshetikhin twist. 
We study the relations between the family
$\{ U_{\cal E(\theta)}(sl(3)) \}$ and the one-dimensional 
set $\{ U_{\cal DJR(\lambda)}(sl(3)) \}$ produced by the standard
Reshetikhin twist from the Drinfeld--Jimbo 
quantization $U_{\cal DJ}(sl(3))$.
These sets of deformations are in one-to-one 
correspondence: each element of $\{ U_{\cal E(\theta)}(sl(3)) \}$ can 
be obtained by a limiting procedure from the unique point in the 
set $\{ U_{\cal DJR(\lambda)}(sl(3)) \}$.  

\end{abstract}
\newpage

%%%%%%%%%%%%%%%%%%%%%%%%%%%%%%%% INTRODUCTION %%%%%%%%%%%%%%%%%%%%%%%%%%%%%%%
\sect{Introduction}

The triangular Hopf algebras and twists (they preserve the triangularity 
\cite{D2,D3}) play an important role in quantum group theory and 
applications \cite{KUL1,VLA,VAL}. Very few types of twists were written 
explicitly in a closed form. The well known example is the jordanian twist 
(${\cal JT}$) of the Borel algebra $B(2)$ ($\{H,E|[H,E]=E\}$) 
with $r=H\wedge E$ \cite{DRIN} where the triangular 
$R$--matrix ${\cal R}=(\faj)_{21}\faj^{-1}$  is defined by the twisting
element \cite{OGIEV,GER}
\begin{equation}
\label{og-twist}
\faj=\exp \{H\otimes \sigma \},
\end{equation}
with $\sigma = \ln (1 + E)$.
In \cite{KLM} it was shown that there 
exist different extensions (${\cal ET}$'s) 
of this twist.
Using the notion of factorizable twist \cite{RSTS} 
the element $\fae \in {\cal U} (sl(N))^{\otimes 2}$,
\begin{equation}
\label{twist-sl(N)}
\fae= \Phi_e \Phi_j =\exp \{2\xi \sum_{i=2}^{N-1}E_{1i}\otimes
E_{iN}e^{-\widetilde{\sigma} }\}\exp \{H\otimes \widetilde{\sigma} \},
\end{equation}
was proved to satisfy the twist equation,
where $E=E_{1N}$, $H=E_{11}-E_{NN}$ is one of the Cartan 
generators $ H \in {\goth h}(sl(N))$, 
$\widetilde{\sigma} =\frac 12\ln (1+2\xi E)$ and
 $\{ E_{ij} \} _{i,j = 1, \dots N} $ is the standard $gl(N)$ basis. 

Studying the family
$\{{\bf L}(\alpha,\beta,\gamma,\delta)_{\alpha + \beta = \delta}\}$
of carrier algebras for extended jordanian 
twists ${\cal F}_{{\cal E}(\alpha,\beta,\gamma,\delta)}$ \cite{V-M}
it is sufficient to consider the one-dimensional set
${\cal L}=\{{\bf L}(\alpha,\beta)_{\alpha + \beta = 1,}\}$  
(for different nonzero $\gamma$'s and $\delta$'s the Hopf 
algebras ${\bf L}_{\cal E}$, obtained by the corresponding twistings, 
are equivalent). 

The connection of the Drinfeld--Jimbo (${\cal DJ}$) 
deformation of a simple Lie algebra $\fra g$ \cite{DRIN,JIMB} with the
jordanian deformation was already pointed out in \cite{GER}. The similarity
transformation of the classical $r$--matrix 
$$
r_{\cal DJ} = \sum_{i=1}^{{\rm rank}(\fra g)}t_{ij}H_i \otimes H_j 
+ \sum_{\alpha \in \Delta_+}
E_{\alpha} \otimes E_{- \alpha}
$$ 
performed by the operator $\exp (v\ {\rm ad}_{E_{1N}})$  turns $r_{\cal DJ}$ 
into  the sum $r_{\cal DJ}+ v\  r_j$ \cite{GER} where 
\begin{eqnarray}
r_j=- v \left( H_{1N}\wedge E_{1N}+2\sum_{k=2}^{N-1}E_{1k}\wedge
E_{kN}\right) .
\end{eqnarray}
Hence, $r_j$ is also a classical $r$--matrix  and  defines the corresponding
deformation. A contraction of the quantum
Manin plane $xy=qyx$ of ${\cal U}_{\cal DJ}(sl(2))$
with the mentioned above similarity transformation in the fundamental
representation $M=1+ v (1-q)^{-1} \rho (E_{12})$ results
in the jordanian plane $x^{\prime }y^{\prime }
=y^{\prime }x^{\prime }+ v {y^{\prime }}^2$ of ${\cal U}_{j}(sl(2))$
\cite{OGIEV}. Thus, the canonical extended 
jordanian twisted algebra 
$U_{{\cal E}(1/2,1/2)}$, which corresponds in our notation to the carrier 
subalgebra ${\bf L}_{(1/2,1/2)}$, can be treated as a limit case for the 
parameterized set of Drinfeld--Jimbo quantizations. Contrary to this fact
other extended twists of $U(sl(N))$ do not reveal such properties with 
respect to the standard deformation. 
In particular, the $U_{\cal P}(sl(4))$ algebra
twisted by the so-called peripheric
twist (${\cal PT}$)  was found to be disconnected with the Drinfeld--Jimbo 
deformation $U_{\cal DJ}(sl(4))$.

In this paper we study the properties of the deformations induced in
$U(sl(3))$ by the set of extended twists \faeab.
We consider the deformations of simple Lie algebras. So, the
parameters $\alpha$ and $\beta$ 
(arising from the reparametrization of the root space)
can be treated as belonging to ${\bf R}^1$. The same is true for other
parameters $(\lambda, \theta, \ldots )$ appearing in this study. 
In the twist equivalence transformations they can be
considered as belonging to ${\bf C}^1$. But in the present approach it is 
sufficient to treat them as real numbers.  

We show that to any Hopf algebra 
$U_{\cal E(\alpha,\beta)}$ one can apply additional
Reshetikhin twist \cite{RES} $\farol$ whose (abelian) carrier 
subalgebra is generated by $K \in {\goth h}(sl(N))$ and $E \in {\bf L}$:
\be
U_{\cal E(\alpha,\beta)} 
\stackrel{{\cal F}_{\widetilde{\cal R}(\lambda)}}{\longrightarrow}
 U_{{\cal E}\widetilde{\cal R}(\alpha,\beta,\lambda)}.
\ee  
However, the carrier subalgebra of ${\cal F}_{\widetilde{\cal
R}(\lambda)}\circ {\cal F}_{{\cal E}(\alpha,\beta)}$
is the same as for ${\cal F}_{{\cal E}(\alpha,\beta)}$  because of the
isomorphism:
\be
   U_{{\cal E}\widetilde{\cal R}(\alpha,\beta,\lambda)}
\approx  U_{\cal E(\alpha +\lambda,\beta -\lambda)}.
\ee
Twists $\farol$ act transitively on
the set $\{ U_{\cal E(\alpha,\beta)} \}$. Simultaneously we consider the
canonical Reshetikhin twist $\fart = 
e^{\theta H_1 \otimes H_2}$ \cite{RES} 
that performs the transition from  $U_{\cal DJ}(sl(3))$  
to the parametric quantization: 
\be
U_{\cal DJ} \stackrel{{\cal F}_{{\cal R}(\theta)}}{\longrightarrow}
 U_{\cal DJR(\theta)}.
\ee
It is worth mentioning that in the case of $U_{\cal DJ}(
gl(3))$ such kind of transformations can be used to obtain possibilities 
for additional twistings  \cite{KM}.
   
Finally, the two sets of parameterized Lie algebras 
are formed: $\{ \gde (\lambda) \}$
and $\{ \gddjr (\theta) \}$. The elements of both of them are dual to $sl(3)$.
Using the technique elaborated in \cite{V-M,KL} we prove a 
one-to-one correspondence between the members of these sets: for 
any $\lambda_0$ fixed there is one and only one $\theta_0$ such 
that $\gde (\lambda_0)$ and $\gddjr (\theta_0)$ are the first order 
deformations of each other. This means that for any $\uneabsl$ there exists 
one and only one such $\undjrtsl$ that these two Hopf algebras can be
connected by a smooth sequence of quantized Lie bialgebras. 

In Section 2 we present a short list of basic relations for twists.
The general properties of extended twists for $U(sl(3))$ are displayed in Section 3. 
There we construct explicitly the peripheric extended twisted algebra
\unpsosl.  In Section 4 the special kind of Reshetikhin
twist for \unpsosl \ is composed and as a result the 
family $\{  U_{{\cal P'}\widetilde{\cal R}(\lambda)}(sl(3))\}$ is 
obtained. We prove that this solves the problem of finding the whole set
$ \{ \uneosl \}$ of extended twists. The relations between 
the multiparametric ${\cal DJ}$ quantizations and twisted 
algebras $ \{\uneosl\}$ are studied in Section 5, and their one-to-one 
correspondence is established.  The defining relations for the canonically
extended twisted algebra $U^{\rm can}_{\cal E}(sl(3))$ are presented in the
Appendix.

%%%%%%%%%%%%%%%%%%%%%%%%%%%%%%%% SECTION 2 %%%%%%%%%%%%%%%%%%%%%%%%%%%%%%%%%

\sect{Basic definitions }

In this section we remind briefly the basic properties of twists.

A Hopf algebra ${\cal A}(m,\Delta ,\epsilon,S)$ with
multiplication $m\colon {\cal A}\otimes {\cal A}\to {\cal A}$,
coproduct $\Delta \colon {\cal A}\to {\cal A}%
\otimes {\cal A}$, counit $\epsilon \colon {\cal A}\to C$,
and antipode $S : {\cal A}\to {\cal A}$ 
can be transformed \cite{D2} by an invertible (twisting) 
element ${\cal F}\in {\cal A} 
\otimes {\cal A}$, ${\cal F}=\sum f_i^{(1)}\otimes f_i^{(2)}$, 
into a twisted
one ${\cal A}_{\cal F}(m,\Delta _{\cal F},\epsilon ,S_{\cal F})$.
This Hopf algebra ${\cal A}_{\cal F}$ has the
same multiplication and counit  but the twisted coproduct and antipode given by
\begin{equation}
\label{def-t}
\Delta _{\cal F}(a)={\cal F}\Delta (a){\cal F}^{-1},\qquad S_{\cal
F}(a)=VS(a)V^{-1},
\end{equation}
with
$$
V=\sum f_i^{(1)}S(f_i^{(2)}), \qquad a\in {\cal A}.
$$
The twisting element has to satisfy the equations
\begin{eqnarray}
\label{def-n}
(\epsilon \otimes  id)({\cal F}) = (id \otimes  \epsilon)({\cal F})=1,
\\[0.2cm]
\label{gentwist}
{\cal F}_{12}(\Delta \otimes  id)({\cal F}) =
{\cal F}_{23}(id \otimes  \Delta)({\cal F}).
\label{TE}
\end{eqnarray}
The first one is just a normalization condition and
follows from the second relation modulo a non-zero scalar factor.

If ${\cal A}$ is a Hopf subalgebra of ${\cal B}$ the twisting 
element ${\cal F}$
satisfying (\ref{def-t})--(\ref{gentwist}) induces the twist
deformation  ${\cal B}_{\cal F}$ of  ${\cal B}$. In this case one can
put $a \in  {\cal B}$ in all the formulas (\ref{def-t}). This will
completely define the Hopf algebra ${\cal B}_{\cal F}$. Let ${\cal A}$ and
${\cal B}$ be the universal enveloping algebras: ${\cal A} = U({\fra l}) 
\subset {\cal B}= U(\g)$ with ${\fra l} \subset \g$. If $U({\fra l})$ 
is the minimal subalgebra on which ${\cal F}$ is completely defined 
as ${\cal F} \in U({\fra l})
\otimes U({\fra l})$ then ${\fra l}$ is called the carrier algebra for 
${\cal F}$ \cite{GER}. 

The composition of appropriate twists can be defined as
${\cal F} = {\cal F}_2 {\cal F}_1$. Here the element ${\cal F}_1$ has to
satisfy the twist equation with the coproduct of the original Hopf algebra,
while ${\cal F}_2$ must be its solution for $\Delta_{{\cal F}_1}$ of the
algebra twisted by ${\cal F}_1$.

If the initial Hopf algebra ${\cal A}$ is quasitriangular with the 
universal element ${\cal R}$ then so is the twisted one 
${\cal A}_{\cal F}(m,\Delta _{\cal F},\epsilon,S_{\cal F},{\cal R}_{\cal F})$
with 
\begin{eqnarray}\label{Rt}
{\cal R}_{\cal F}=  {\cal F}_{21} \,{\cal R} \,{\cal F}^{-1}.
\end{eqnarray}

Most of the explicitly known twisting elements have the factorization
property with respect to comultiplication
$$
(\Delta \otimes id)({\cal F})={\cal F}_{23}{\cal F}_{13}\qquad \mbox{or}
\qquad
(\Delta \otimes id)({\cal F})={\cal F}_{13}{\cal F}_{23}\,,
$$
and
$$
(id \otimes \Delta)({\cal F})={\cal F}_{12}{\cal F}_{13}\qquad \mbox{or}
\qquad
(id \otimes \Delta)({\cal F})={\cal F}_{13}{\cal F}_{12}\,.
$$
To guarantee the validity of the twist equation, these identities are to be 
combined with  the additional requirement
${\cal F}_{12}{\cal F}_{23}={\cal F}_{23}{\cal F}_{12}$
or the Yang--Baxter equation on ${\cal F}$ \cite{RSTS}.

An important subclass of factorizable twists consists of elements
satisfying the equations
\begin{eqnarray} \label{f-twist1}
(\Delta \otimes id)({\cal F})={\cal F}_{13}{\cal F}_{23}\,,   
\\ [0,2cm] \label{f-twist2}
(id\otimes \Delta _{\cal F})({\cal F})={\cal F}_{12}{\cal F}_{13 }\,.
\end{eqnarray}
Apart from the universal $R$--matrix ${\cal R}$ that satisfies these
equations for $\Delta_{\cal F}=\Delta ^{op}$ ($\Delta ^{op}=\tau\circ \Delta$, 
where $\tau(a\otimes b)=b\otimes a$)  there are two
more well developed  cases of such twists: the jordanian twist of  a Borel 
algebra $B(2)$   where ${\cal F}_j$
has the form (\ref{og-twist}) (see \cite{OGIEV}) and the extended
jordanian  twists (see \cite{KLM} and \cite{V-M} for details).

According to the result by Drinfeld \cite{D3}  skew (constant) 
solutions of the classical Yang--Baxter equation (CYBE) can be quantized and 
the  deformed algebras thus  obtained can be presented in a form of twisted 
universal enveloping  algebras. On the other hand, such solutions of CYBE 
can be connected with  the quasi-Frobenius carrier subalgebras of the initial 
classical Lie  algebra \cite{STO}. A Lie  algebra $\g(\mu)$, 
with the Lie composition $\mu$, 
is called Frobenius if there exists a linear functional $g^* \in \g^*$
such that the form $b(g_1,g_2)=g^*(\mu(g_1,g_2))$ is nondegenerate.
This means that $\g$ must have a nondegenerate 2--coboundary $b(g_1,g_2) \in
B^2(\g,{\bf K})$. The algebra is called quasi-Frobenius if it has a 
nondegenerate 2--cocycle $b(g_1,g_2) \in Z^2(\g,{\bf K})$ (not 
necessarily a coboundary). The classification of quasi-Frobenius
subalgebras in $sl(n)$ was given in \cite{STO}.

The deformations of quantized algebras include the deformations of their
Lie bialgebras $(\g,\g^*)$. Let us fix one of the constituents
$\g^*_1(\mu^*_1)$ (with composition $\mu^*_1$) and deform it  in the  first
order
$$
(\mu^*_1)_t = \mu^*_1 + t\mu^*_2,
$$
its deforming function $\mu^*_2$ is also a Lie product and the deformation
property becomes reciprocal: $\mu^*_1$ can be considered as a first order 
deforming function  for algebra $\g^*_2(\mu^*_2)$. Let $\g(\mu)$
be a Lie algebra that form Lie bialgebras with  both
$\g^*_1$ and $\g^*_2$. This means that we have a one-dimensional 
family $\{ (\g,(\g^*_1)_t) \}$ of Lie bialgebras and correspondingly a one
dimensional family of quantum deformations $ \{{\cal A}_t(\g,(\g^*_1)_t) \}$ 
\cite{ETI}. This situation provides the possibility to construct in the 
set of Hopf algebras a smooth curve connecting quantizations of the 
type ${\cal A}(\g,\g^*_1)$ with those of ${\cal A}(\g,\g^*_2)$. Such  smooth 
transitions can involve contractions provided
$\mu^*_2 \in B^2 (\g^*_1,\g^*_1)$. This happens in the case of
${\cal JT,\ ET}$ and some other twists (see \cite{KL} and references therein).

%%%%%%%%%%%%%%%%%%%%%%%% SECTION 3 %%%%%%%%%%%%%%%%%%%%%%%%%%%%%%%%%%%%%% 
\sect{Extended twist for $U(sl(3))$}

Extended jordanian twists are associated with the set 
$\{{\bf L}(\alpha,\beta,\gamma,\delta)_{\alpha + \beta = \delta}\}$
of Frobenius algebras \cite{KLM},\cite{V-M}
\begin{equation}
\begin{array}{l}
 [H,E] = \delta E, \quad [H,A] = \alpha A, \quad [H,B] = \beta B, \\[0.2cm]
[A,B] = \gamma E, \quad[E,A] = [E,B] = 0, \quad \quad
\alpha + \beta = \delta .
\end{array}
\end{equation}
For limit values of $\gamma$ and $\delta$ the structure of ${\bf L}$ 
degenerates. For the internal (nonzero) values of $\gamma$ and $\delta$ the 
twists associated with the corresponding ${\bf L}$'s are equivalent. It is 
sufficient to study the one-dimensional subvariety 
${\cal L}=\{{\bf L}(\alpha,\beta)_{\alpha + \beta = 1}\}$, 
that is to consider the carrier algebras 
\begin{equation}
\begin{array}{l}
 [H,E] = E, \quad [H,A] = \alpha A, \quad [H,B] = \beta B, \\[0.2cm]
[A,B] = E, \quad[E,A] = [E,B] = 0, \quad \quad
\alpha + \beta = 1 .
\end{array}
\end{equation}
The corresponding group 2--cocycles (twists) are
\be
\label{t-ext}
\faeab =  \Phi_{\cal E(\alpha, \beta)} \Phi_j
\ee  
or 
\be
\label{t-ext-s}
\faesab = \Phi_{\cal E'(\alpha, \beta)}  \Phi_j
\ee  
with
\begin{equation}
\label{fractions}
\begin{array}{lcr}
 \Phi_j & = & \faj = \exp \{H\otimes \sigma \}, \\[0.2cm]
 \Phi_{\cal E(\alpha, \beta)} & = & 
\exp \{ A \otimes B e^{-\beta \sigma} \}, \\[0.2cm]
 \Phi_{\cal E'(\alpha, \beta)} & = & 
\exp \{ -B \otimes A e^{-\alpha \sigma} \}.
\end{array}
\end{equation}
Twists (\ref{t-ext}) and (\ref{t-ext-s})
 define the deformed Hopf algebras ${\bf L}_{\cal E(\alpha, \beta)} $ 
with the co-structure
\begin{equation}
\label{e-costr}
\begin{array}{lcl}
 \deleabh & = & H \otimes e^{-\sigma} + 1 \otimes H 
              - A \otimes B e^{-(\beta + 1)\sigma}, \\[0.2cm]
 \deleaba & = & A \otimes e^{-\beta \sigma} + 1 \otimes A , \\[0.2cm]
 \deleabb & = & B \otimes e^{\beta \sigma} + e^{\sigma} \otimes B , \\[0.2cm]
 \deleabe & = & E \otimes e^{\sigma} + 1 \otimes E ;
\end{array}
\end{equation}
 and   
$ {\bf L}_{\cal E'(\alpha, \beta)} $ defined by 
\begin{equation}
\label{e-pr-costr}
\begin{array}{lcl}
 \delesabh & = & H \otimes e^{-\sigma} + 1 \otimes H 
              + B \otimes A e^{-(\alpha + 1)\sigma}, \\[0.2cm]
 \delesaba & = & A \otimes e^{\alpha \sigma} + e^{\sigma} \otimes A , 
\\[0.2cm]
 \delesabb & = & B \otimes e^{-\alpha \sigma} + 1 \otimes B , \\[0.2cm]
 \delesabe & = & E \otimes e^{\sigma} + 1 \otimes E.
\end{array}
\end{equation}
The sets  $\{ {\bf L}_{\cal E(\alpha, \beta)} \}$ and
$\{ {\bf L}_{\cal E'(\alpha, \beta)} \}$
 are equivalent due to the Hopf isomorphism 
${\bf L}_{\cal E(\alpha, \beta)} \approx {\bf L}_{\cal E'(\beta, \alpha)}$:
\be
\{ {\bf L}_{\cal E}(\alpha,\beta) \} \approx \{ {\bf L}_{\cal E'}
(\alpha,\beta) \} 
\approx \{ {\bf L}_{\cal E}(\alpha \geq \beta) \} \cup \{ {\bf L}_{\cal E'}
(\alpha \geq \beta) \}.
\ee 
So, it is sufficient to use only one of the extensions 
either $\Phi_{\cal E(\alpha, \beta)}$ or $\Phi_{\cal E'(\alpha, \beta)}$,
or a half of the domain for $(\alpha, \beta)$.

The set ${\cal L}=\{{\bf L}(\alpha,\beta)_{\alpha + \beta = 1,}\}$  is just 
the family of 4-dimensional Frobenius algebras that one finds in $U(sl(3))$
\cite{STO}. It was  mentioned in \cite{KLM} that complicated 
calculations are needed to write down all the defining coproducts for the
canonical extended twisted $U_{\cal E}^{\rm can}(sl(3))$. Here we shall show 
how to overcome partially this difficulty and to get all the defining 
relations in the explicit form. 

First we shall construct the simplest member of the family $\{ \uneabsl \}$
--- one of the peripheric twisted algebras \unpsosl. Then, the 
additional parameterized twist  will be applied and finally 
we shall prove that the whole set 
$\{ \uneabsl \}$ is thus obtained.

Consider the subalgebra $\L(0,1) \subset sl(3)$ with generators
\be
\begin{array}{l}
H  =  \frac{1}{3} (H_{13}+ H_{23}) = \frac{1}{3}(E_{11} 
+ E_{22} -2E_{33}),\\[0.2cm]
A  =   E_{12}, \quad  B  =   E_{23}, \quad  E  =   E_{13},
\end{array}
\ee
and the compositions  
\be
\label{perl}
\begin{array}{l}
[H, E_{13} ]   =  E_{13},\quad [H, E_{12} ]   =  0, 
\quad [H, E_{23} ]   =  E_{23}, \\[0.2cm]
[E_{12}, E_{23} ]   =  E_{13}, \quad [E_{12}, E_{13} ]   
=  [E_{23}, E_{13} ] = 0 .
\end{array}
\ee
According to the results obtained in \cite{V-M} (see  
formulas (\ref{t-ext-s}) and (\ref{fractions}))
one of the peripheric twists attributed to this algebra has the form
\be
\faps = \Phi_{\cal P'}\Phi_j = e^{- E_{23} \otimes E_{12}} 
e^{H \otimes \sigma}.
\ee 
Applying to $U(sl(3))$ the twisting procedure with \faps we construct 
the Hopf algebra \unpsosl \ with  
the usual multiplication  of $U(sl(3))$ and the coproduct defined by the 
relations:
\be
\label{up-co-m}
\begin{array}{lcl}
\delps(H_{12})&=&  H_{12} \otimes 1 + 1 \otimes H_{12} + H \otimes 
(\exms - 1) +  E_{23} \otimes E_{12}\exms,\\[2mm]
\delps(H_{13})&=& H_{13} \otimes 1 + 1 \otimes H_{13} + 2H \otimes (\exms -1) 
 + 2 E_{23} \otimes E_{12}\exms, \\[2mm]                        
\delps(E_{12})&=& E_{12} \otimes 1 + \exs \otimes E_{12},\\[2mm]
\delps(E_{13})&=& E_{13} \otimes \exs + 1 \otimes E_{13},\\[2mm]
\delps(E_{21})&=& E_{21} \otimes 1 + 1 \otimes E_{21} - H \otimes 
E_{23} \exms - E_{23} \otimes H_{12} \\[1mm]
             &  & \ \ - \ E_{23} \otimes E_{12}E_{23}\exms  
              +  H E_{23} \otimes (1- \exms) 
              - E_{23}^2 \otimes E_{12} \exms, \\[2mm] 
\delps(E_{23})&=& E_{23} \otimes 1 + 1 \otimes E_{23},\\[2mm]
\delps(E_{31})&=& E_{31} \otimes \exms + 1 \otimes E_{31} 
              + H \otimes H_{13} \\[1mm] 
             & &   \ \ + \  (1 - H) H \otimes (\exms - \exmds)\\[1mm]
													& & \ \ + \  (1 -  H) E_{23} \otimes E_{12}(\exms - 2\exmds) 
 -  E_{21} \otimes E_{12} \exms \\[1mm]
             & & \ \  + \  E_{23} \otimes E_{32}  
              +  E_{23} \otimes H_{13} E_{12} \exms
              +   E_{23}^2 \otimes E_{12}^2 \exmds, \\[2mm]

\delps(E_{32})&=& E_{32} \otimes \exms + 1 \otimes E_{32} 
             + (H - H_{23}) \otimes E_{12} \exms. \\[1mm] 
\end{array}
\ee

The universal ${\cal R}$--matrix for this algebra is
\be
{\cal R}_{\cal P'} = e^{- E_{12} \otimes E_{23}} e^{\sigma \otimes H} 
 e^{- H \otimes \sigma}  e^{ E_{23} \otimes E_{12}},
\ee 
and the classical $r$--matrix can be written in the form
\be
\label{rmat}
r_{\cal P'} =  E_{23} \wedge E_{12} + \frac{1}{3}E_{13} 
            \wedge (E_{11} + E_{22} - 2E_{33}).
\ee

By means of this $r$--matrix (or directly from the
coproducts (\ref{up-co-m})) the following Lie compositions for 
\gdps (the algebra dual to $sl(3)$ in this quantization) can be obtained
\be
\label{p-dual}
\begin{array}{lclrcl}
[X_{11}, X_{13} ]  & = & -\frac{1}{3}(X_{11}-X_{33}), &\qquad
[X_{12}, X_{23} ]  & = & - (X_{11}-X_{33}),\\[1mm]
[X_{22}, X_{13} ]  & = & -\frac{1}{3}(X_{11}-X_{33}), &
[X_{12}, X_{13} ]  & = & - X_{12},\\[1mm]
[X_{33}, X_{13} ]  & = & +\frac{2}{3}(X_{11}-X_{33}), &
[X_{12}, X_{21} ]  & = &  X_{31},\\[1mm]
[X_{11}, X_{23} ]  & = & +\frac{2}{3} X_{21}, &
[X_{13}, X_{31} ]  & = &  X_{31},\\[1mm]
[X_{22}, X_{23} ]  & = & -\frac{4}{3} X_{21}, &
[X_{23}, X_{32} ]  & = &  X_{31},\\[1mm]
[X_{33}, X_{23} ]  & = & +\frac{2}{3} X_{21}, &
[X_{13}, X_{32} ]  & = & + X_{32},\\[1mm]
[X_{11}, X_{33} ]  & = & +\frac{1}{3} X_{31}, &
[X_{22}, X_{33} ]  & = & -\frac{1}{3} X_{31},\\[1mm]
[X_{11}, X_{12} ]  & = & +\frac{1}{3} X_{32}, &
[X_{12}, X_{33} ]  & = & -\frac{1}{3} X_{32},\\[1mm]
[X_{11}, X_{22} ]  & = & -\frac{1}{3} X_{31}, &
[X_{12}, X_{22} ]  & = & +\frac{2}{3} X_{32}.\\[1mm]
\end{array}
\ee

%%%%%%%%%%%%%%%%% SECTION 4 %%%%%%%%%%%%%%%%%%%%%%%%%%%%%%%%%%%%%%%
\sect{Reshetikhin twist action on $U_{\cal E}(sl(3))$}

The main observation with respect to our present aim is
that besides the primitive element $\sigma$ the twisted algebra \unpsosl \
contains the primitive Cartan generator $K$ 
\be
K = \frac{1}{3}(H_{12} - H_{23}).
\ee
The element $K^*$ dual to $K$ is orthogonal to the root $E^*$ of $E 
\in \L(\alpha,\beta)$, that is, $K$ commutes with $\sigma$.
So  \unpsosl \ contains the Abelian subalgebra
\begin{equation}
\begin{array}{lcl}  
\delps(K) & = & K \otimes 1 + 1 \otimes K, \\[2mm]
\delps(\sigma) & = & \sigma \otimes 1 + 1 \otimes \sigma. 
\end{array}
\quad [K, \sigma]  =  0,
\end{equation}   
Thus, the additional Reshetikhin twist 
\be
\farol = e^{\lambda K \otimes \sigma}
\ee
is applicable to the previously obtained Hopf algebra,
\be
\unpsosl \stackrel{\farol}{\longrightarrow}
 U_{{\cal P'}\widetilde{\cal R}(\lambda)}(sl(3)).
\ee
The new twisted algebra $ U_{{\cal P'}\widetilde{\cal R}(\lambda)}(sl(3))$ is 
defined by the relations:
 \be
\label{ug-co-m}
\begin{array}{lcl}
\delpsol (H_{12})&=&  H_{12} \otimes 1 + 1 \otimes H_{12} 
                  + (\lambda K + H) \otimes (\exms - 1) \\[1mm]
             & & \ \ + \ E_{23} \otimes E_{12}e^{-(\lambda + 1) \sigma},\\[2mm]
\delpsol (H_{13})&=& (H_{13} - 2(\lambda K + H)) \otimes 1 
                 +2 (\lambda K + H) \otimes \exms \\[1mm]
             & & \ \  +1\ \otimes H_{13}   
              +2 E_{23} \otimes E_{12} e^{-(\lambda +1)\sigma}, \\[2mm]
\delpsol(E_{12})&=& E_{12} \otimes e^{\lambda \sigma} 
                 + \exs \otimes E_{12},\\[2mm]
\delpsol(E_{13})&=& E_{13} \otimes \exs + 1 \otimes E_{13},\\[2mm]
\delpsol(E_{21})&=& E_{21} \otimes  e^{-\lambda \sigma} + 1 \otimes E_{21} \\[1mm]
             & & \ \ - \ E_{23} \otimes H_{12} e^{-\lambda \sigma}
                 - (\lambda K + H) \otimes E_{23} \exms   \\[1mm]
             & &\ \ +\ (\lambda K + H) E_{23} \otimes 
                 (e^{-\lambda \sigma} -  e^{-(\lambda +1)\sigma})  \\[1mm]
             & & \ \ - \ E_{23}^2 \otimes E_{12} e^{-(2\lambda +1)\sigma}
               - E_{23} \otimes E_{12}E_{23} e^{-(\lambda +1)\sigma},  \\[2mm] 
\delpsol(E_{23})&=& E_{23} \otimes e^{- \lambda \sigma} + 1 \otimes E_{23} ,\\[2mm]
\delpsol(E_{31})&=& E_{31} \otimes \exms + 1 \otimes E_{31} 
              + (\lambda K + H) \otimes H_{13} \exms  \\[1mm]
             & & \ \ +\  E_{23} \otimes E_{32} e^{- \lambda \sigma} \\[1mm] 
             & & \ \ + \ (1 - \lambda K - H)(\lambda K + H) \otimes (\exms - \exmds)
\\[1mm]
             & &\ \  - \ E_{21} \otimes E_{12} e^{-(\lambda +1)\sigma} \\[1mm] 
             & &\ \  + \ (1 - \lambda K - H)E_{23} \otimes E_{12}e^{-\lambda \sigma}
                 (\exms - 2\exmds) \\[1mm]  
             & & \ \ + \ E_{23} \otimes H_{13} E_{12}e^{-(\lambda +1)\sigma} 
                 + E_{23}^2 \otimes E_{12}^2 e^{-2(\lambda +1)\sigma} , \\[2mm]
\delpsol(E_{32})&=& E_{32} \otimes  e^{(\lambda -1)\sigma} + 1 \otimes E_{32} 
                 + (\lambda + 1) K \otimes E_{12} \exms. \\[1mm] 
\end{array}
\ee

According to the associativity of twisting transformations 
the same parameterized set of algebras could be obtained directly 
from $U(sl(3))$ using the composite twist
\be
\fapsrol = \farol \Phi_{\cal E'} \Phi_{j} = e^{\lambda K \otimes \sigma}
 e^{- E_{23} \otimes E_{12}} e^{H \otimes \sigma}.
\ee
This twisting element can be written in the form
\be\label{twistpr}
\fapsrol  = 
 e^{- E_{23} \otimes E_{12}e^{-\lambda \sigma}} e^{(H + \lambda K) \otimes 
\sigma}.
\ee
The latter is the extended twist for the Lie algebra
\begin{equation}\label{liealgebra}
\begin{array}{l}
 [H + \lambda K,E_{13}] = E_{13}, \quad 
 [H + \lambda K,E_{12}] = \lambda E_{12}, \quad 
 [H + \lambda K,E_{23}] = (1 - \lambda) E_{23}, \\[0.2cm]
[E_{12},E_{23}] = E_{13}, \quad[E_{12},E_{13}] = [E_{23},E_{13}] = 0 .
\end{array}
\end{equation}
The relations (\ref{twistpr}) and (\ref{liealgebra}) signify that the family 
$\{ U_{{\cal P'}\widetilde{\cal R}(\lambda)} (sl(3)) \}$ is the complete 
set of
twisted Hopf algebras related to the Frobenius 
subalgebras $\{ {\bf L}_{\cal E(\alpha, \beta)} 
             \in sl(3) \}$ and that $\lambda = \alpha $.

It must be also stressed that the appropriate Reshetikhin twist of the 
type $\farol$ can be constructed for any algebra  
$U_{{\cal E(\alpha, \beta)}}(sl(3))$ --- there always exists a Cartan 
element whose dual is orthogonal to the root $\nu_{E}$.

Note that any triple of roots $\{ \alpha, \beta, \gamma \mid \alpha + \beta =
 \gamma \}$ of  the $sl(3)$ root system  can play the role of the triple 
$\{ \nu_{12}, \nu_{23}, \nu_{13} \}$ that was selected in our case to form 
the carrier subalgebra. The formulas above 
are irrelevant to this choice, only the interrelations of roots
are important. In $sl(3)$ there always exists the equivalence transformation of 
the root system that identify any such triple with the fixed one.

The obtained set of Hopf algebras corresponds to the parameterized family
$r_{\cal E'(\theta)}$ of $r$--matrices 
\be
\label{est-rmat}
r_{\cal E'(\theta)}
       =  E_{23} \wedge E_{12} + \frac{1}{2}E_{13} \wedge H_{13} +
            \frac{1}{2}\theta E_{13} \wedge (H_{12} - H_{23}),
\ee
where we use the  parameter $\theta = \frac{1}{3}(2\lambda - 1)$ measuring the
deviation of the extended twist from the canonical rather than from the
peripheric one.

Algebras $ U_{{\cal P'}\widetilde{\cal R}(\lambda)}(sl(3))$ are the 
quantizations
of the Lie bialgebras $(sl(3),{\fra g}_{\cal E'(\theta)}^*)$. The compositions of
\gdest are  easily derived with the help of (\ref{est-rmat}): 
\be
\label{es-dual}
\begin{array}{lclrcl}
[X_{11}, X_{12} ]  & = & \frac{1}{2}(1 + \theta) X_{32},
&\qquad [X_{11}, X_{22} ]  & = & \theta X_{31}, \\[1mm]
[X_{11}, X_{23} ]  & = & \frac{1}{2}(1 - \theta)X_{21},
&[X_{11}, X_{13} ]  & = & - \frac{1}{2}(\theta + 1)(X_{11} - X_{33}),\\[1mm]
[X_{11}, X_{33} ]  & = &  - \theta X_{31}, 
&[X_{12}, X_{13} ]  & = & \frac{1}{2}(3\theta - 1) X_{12},\\[1mm]
[X_{12}, X_{21} ]  & = &  X_{31},
&[X_{12}, X_{23} ]  & = & - (X_{11}-X_{33}),\\[1mm]
[X_{12}, X_{22} ]  & = &  (\theta + 1) X_{32},
&[X_{12}, X_{33} ]  & = & - \frac{1}{2}(\theta + 1) X_{32}, \\[1mm]
[X_{13}, X_{21} ]  & = &  \frac{1}{2}(3\theta + 1) X_{21},
&[X_{13}, X_{22} ]  & = &  -\theta(X_{11} - X_{33}),\\[1mm]
[X_{13}, X_{23} ]  & = &  \frac{1}{2}(3\theta + 1) X_{23},
&[X_{13}, X_{31} ]  & = & X_{31},\\[1mm]
[X_{13}, X_{32} ]  & = &  \frac{1}{2}(1 - 3\theta) X_{32},
&[X_{13}, X_{33} ]  & = & \frac{1}{2}(\theta - 1) (X_{11} - X_{33}),\\[1mm]
[X_{22}, X_{23} ]  & = &  (\theta - 1) X_{21},
&[X_{22}, X_{33} ]  & = &  \theta X_{31},\\[1mm]
[X_{23}, X_{32} ]  & = &  X_{31},
&[X_{23}, X_{33} ]  & = & \frac{1}{2}(\theta - 1) X_{21}.\\[1mm]
\end{array}
\ee

%%%%%%%%%%%%%%%%% SECTION 5 %%%%%%%%%%%%%%%%%%%%%%%%%%%%%%%%%%%%%%%
\sect{Multiparametric Drinfeld--Jimbo and ${\cal ET}$ quantizations}

The twisting element for the Reshetikhin twist \cite{RES} for \undjosl, 

\be
\far = e^{\eta H_{23} \wedge H_{12}},
\ee
converts \undjosl \ into the twisted algebra \undjrosl with the 
$r$--matrix of the form
\be 
r_{\cal DJR} = r_{\cal DJ} + r_{\cal R} = r_{\cal DJ} + \eta H_{12} 
\wedge H_{23} = 
 r_{\cal DJ} + \eta (E_{11}\wedge E_{33} - E_{11}\wedge E_{22} - 
E_{22}\wedge E_{33}).
\ee  
This signifies that the corresponding dual Lie algebra
\gddjr is the first order deformation of \gddj by \gdr and 
$\eta$ can be viewed  as a deformation parameter. The nonzero compositions 
of \gdr are the following ones :
\be
\label{r-dual}
\begin{array}{lclrcl}
[X_{11}, X_{12} ]  & = &- X_{12},
&\qquad [X_{11}, X_{21} ]  & = &  X_{21}, \\[1mm]
[X_{22}, X_{12} ]  & = & -X_{12},
&[X_{22}, X_{21} ]  & = &  X_{21},\\[1mm]
[X_{33}, X_{12} ]  & = & 2 X_{12}, 
&[X_{33}, X_{21} ]  & = & -2 X_{21},\\[1mm]
[X_{11}, X_{13} ]  & = &  X_{13},
&[X_{11}, X_{31} ]  & = & -X_{31},\\[1mm]
[X_{22}, X_{13} ]  & = & -2 X_{13},
&[X_{22}, X_{31} ]  & = & 2 X_{31},\\[1mm]
[X_{33}, X_{13} ]  & = &   X_{13},
&[X_{33}, X_{31} ]  & = & -X_{31},\\[1mm]
[X_{11}, X_{23} ]  & = & 2 X_{23},
&[X_{11}, X_{32} ]  & = & -2 X_{32},\\[1mm]
[X_{22}, X_{23} ]  & = & -X_{23},
&[X_{22}, X_{32} ]  & = &  X_{32},\\[1mm]
[X_{33}, X_{23} ]  & = & -X_{23},
&[X_{33}, X_{32} ]  & = &  X_{32}.\\[1mm]
\end{array}
\ee 

The compositions \mdjr \, of the   
algebra  \gddjr that was deformed in the first order by \mre \, are:
\be
\label{djr-dual-3}
\begin{array}{lclrcl}
[X_{11}, X_{12} ]  & = & X_{12} - \eta X_{12},
&\qquad [X_{11}, X_{21} ]  & = & X_{21} + \eta X_{21},\\[1mm]
[X_{11}, X_{13} ]  & = & X_{13} + \eta X_{13},
&[X_{11}, X_{31} ]  & = & X_{31} - \eta X_{31},\\[1mm]
[X_{11}, X_{23} ]  & = &  2\eta X_{23}, 
&[X_{11}, X_{32} ]  & = &  -2\eta X_{32},\\[1mm]
[X_{22}, X_{12} ]  & = & - X_{12} - \eta X_{12},
&[X_{22}, X_{21} ]  & = & - X_{21} + \eta X_{21},\\[1mm]
[X_{22}, X_{13} ]  & = &  -2 \eta X_{13},
&[X_{22}, X_{31} ]  & = &  2 \eta X_{31},\\[1mm]
[X_{22}, X_{23} ]  & = &  X_{23} - \eta X_{23},
&[X_{22}, X_{32} ]  & = &  X_{32} + \eta X_{32},\\[1mm]
[X_{33}, X_{12} ]  & = & 2 \eta X_{12},
&[X_{33}, X_{21} ]  & = & - 2 \eta X_{21},\\[1mm]
[X_{33}, X_{13} ]  & = & - X_{13} + \eta X_{13},
&[X_{33}, X_{31} ]  & = & - X_{31} - \eta X_{31},\\[1mm]
[X_{33}, X_{23} ]  & = & - X_{23} - \eta X_{23},
&[X_{33}, X_{32} ]  & = & - X_{32} + \eta X_{23},\\[1mm]
[X_{12}, X_{23} ]  & = & 2 X_{13},
&[X_{21}, X_{32} ]  & = & 2 X_{31}.\\[1mm]
\end{array}
\ee

According to the lemma proved in  \cite{V-M} the necessary and sufficient 
condition for the existence of a smooth transition connecting two quantized
Lie bialgebras $U_q(\g,\g^*_1)$ and $U_q(\g,\g^*_2)$ is the existence of the 
first
order deformation of $\mu^*_1$ by $\mu^*_2$ (and vice versa). In our case 
this is the combination of compositions (\ref{es-dual}) 
and (\ref{djr-dual-3}),
\be
\mu^*(s,t) = s \mdjr (\eta) + t \mes (\theta),
\ee 
that must be checked. The direct computations show that $\mu^*(s,t)$ is 
a Lie composition if and only if $\eta = \theta$.

Thus we have proved that for any \uneabsl\, there exists one and only one 
twisted Drinfeld--Jimbo deformation $\undjrlsl$ that can be connected with
the twisted algebra by a smooth path whose points are the deformation
quantizations.

Remember that both $\mdjr (\eta)$ and $\mes (\theta)$ are the linear
combinations of Lie compositions ($\mdj$ and $\mre$, $\mps$ and $\mret$). So,
we have a four-dimensional space of compositions with two fixed planes of Lie
compositions containing correspondingly $\mdjr (\eta)$ and $\mes (\theta)$.
From these two planes only the correlated lines (with $\eta = \theta$) 
belong to the Lie subspaces that intersect both planes.

%%%%%%%%%%%%%%%%%     SECTION CONCLUSIONS      %%%%%%%%%%%%%%%%%%%%%%%%%%%%%
\sect{Conclusions}

The $r$--matrix $r_{\cal DJR}(\eta)$ can be transformed into 
the mixed $r$--matrix $r_{\cal DJR}(\eta) + v \ r_{\cal E(\eta)}$ 
with the help of an operator $\exp\{ v\ {\rm ad}_{E} \}$ similarly to
the  ordinary case when $r_{\cal DJ}$ is transformed 
into $r_{\cal DJ} + v \ r_{\cal E}^{\rm can}$ \cite{GER}. We want to note
that the element $E$ may correspond to any root $\nu$ of the
$sl(3)$ root system. Varying the roots one shall arrive at the
$r$--matrices attributed  to different (though equivalent) sets of
extended twisted algebras.

The canonically extended twisted  algebra $U^{\rm can}_{\cal E}(sl(3))$ 
introduced in \cite{KLM} is a special case of extended 
twisted algebras $\{ U_{{\cal E'(\alpha, \beta)}}(sl(3)) \}$. It corresponds 
to the situation when the functional $H^*$ is parallel to the root $\nu_{E}$. 
For the Lie algebras of $A_n$ series this means that $\alpha = \beta = 1/2$. 
In the Appendix we present the full table of the defining relations for this
Hopf algebra.       

The peripheric twists helped us to obtain the explicit form of 
the comultiplication
for all the extended twisted Hopf algebras originated from
$U(sl(3))$. In the set $\{ \uneabsl\}$ algebras produced by peripheric twists 
were not distinguished by 
their relations neither with Drinfeld--Jimbo twists ($\{ \undjosl\}$) nor
with Reshetikhin twists. We want to note that the situation changes when
one studies the specific properties of extensions for peripheric twisted
algebras.

The construction presented in this paper can be performed for any 
two-dimensional sublattice of the root lattice of any simple Lie algebra.
For any highest root of the ``triple" there exists the Cartan generator
whose dual is orthogonal to this root. This means that the corresponding
special Reshetikhin twist can always be constructed. The same is true 
also for the so called
special injections of ${\cal L} \in \g $. In this case the ``triple" will
be realized in the root space submerged in that of the initial simple
algebra. Whatever the injection is an ordinary Reshetikhin twist can be 
applied to the $U_{\cal DJ}(\g)$ to coordinate the properties 
of $U_{\cal DJR}(\g)$ and $U_{\cal ER}(\g)$. 
The extended twists for different injections and the role of the peripheric 
twists will be studied in detail in a forthcoming publication.   

%%%%%%%%%%%%%%%%%%%%%%%%%%%%%%%%%%%%%%%%%%%%%%%%%%%%%%%%%%%%%%%%%%%%%%
\section*{Acknowledgments}   

The authors are thankful to Prof. P.P.Kulish for his important comments.
V. L. would like to thank the DGICYT of the Ministerio de Educaci\'on y 
Cultura de  Espa\~na for supporting his sabbatical stay (grant
SAB1995-0610). This  work has been partially supported by DGES of the 
Ministerio de  Educaci\'on y Cultura of  Espa\~na under Project
PB95-0719, the Junta de Castilla y Le\'on (Espa\~na) and the  Russian 
Foundation for Fundamental Research under grant  97-01-01152.

%%%%%%%%%% the bibliography %%%%%%%%%%%%%%%%%%%%%%%%%%%%%%%%%%%%%%%%%%

%%%%%%%%%%%%%%%%%     APPENDIX     %%%%%%%%%%%%%%%%%%%%%%%%%%%%%%
\section*{Appendix}  

In \cite{KLM} the ${\cal E}$-twisted algebra $U^{\rm can}_{\cal E}(sl(3))$ 
 was introduced and some of its comultiplications
where presented explicitly. In the family 
$\{ U_{{\cal E'(\alpha, \beta)}}(sl(3)) \}$ it corresponds to the case 
$\alpha = 1/2$. The involution
\be
\label{inv-}
\begin{array}{lcllcl}
 E_{12} & \rightleftharpoons & E_{23}, \quad &  E_{32} & 
\rightleftharpoons & - E_{21},\\
 E_{23} & \rightleftharpoons & - E_{12}, \quad & H_{12} & 
\rightleftharpoons & H_{23},\\
 E_{21} & \rightleftharpoons & E_{32}, \quad &  H_{23} & 
\rightleftharpoons & H_{12},
\end{array}
\ee
transforms $ U_{{\cal E'}(1/2, 1/2)}(sl(3)) $ 
into $ U_{{\cal E}(1/2, 1/2)}(sl(3)) $.
The full list of defining coproducts 
for  $U^{\rm can}_{\cal E}(sl(3))$ can be thus obtained:
\be
\label{ue-co-m}
\begin{array}{lcl}
\dele^{\rm can}(H_{23})&=&  H_{23} \otimes 1 + 1 \otimes H_{23} 
                  + \frac{1}{2}H_{13} \otimes (\exmdst - 1) 
                  - 2 \xi E_{12} \otimes E_{23}e^{-3 \widetilde{\sigma}}
,\\[2mm]
\dele^{\rm can}(H_{13})&=&  H_{13} \otimes \exmdst
                 +1 \otimes H_{13} 
                 -4 \xi E_{12} \otimes E_{23} e^{-3 \widetilde{\sigma}}, 
\\[2mm]
\dele^{\rm can}(E_{23})&=& E_{23} \otimes e^{\widetilde{\sigma}} 
                 + \exdst \otimes E_{23},\\[2mm]
\dele^{\rm can}(E_{13})&=& E_{13} \otimes \exdst + 1 \otimes E_{13},
\\[2mm]
\dele^{\rm can}(E_{32})&=& E_{32} \otimes  e^{- \widetilde{\sigma}} 
+ 1 \otimes E_{32}
                 +2\xi E_{12} \otimes H_{23} e^{- \widetilde{\sigma}}
                 +\xi H_{13} \otimes E_{12} \exmdst  \\[1mm]
             & &\ \  -\ \xi H_{13} E_{12} \otimes 
                 (e^{-\widetilde{\sigma}} -  e^{-3\widetilde{\sigma}}) 
                 - 4 \xi^2 E_{12}^2 \otimes E_{23} e^{-4 \widetilde{\sigma}}
                 \\[1mm]
             & &\ \  -\ 4 \xi^2 E_{12} \otimes E_{23}E_{12} 
e^{-3 \widetilde{\sigma}},   \\[2mm] 
\dele^{\rm can}(E_{12})&=& E_{12} \otimes e^{- \widetilde{\sigma}} 
+ 1 \otimes E_{12},\\[2mm]
\dele^{\rm can}(E_{31})&=& E_{31} \otimes \exmdst + 1 \otimes E_{31} 
              + \xi  H_{13} \otimes H_{13} \exmdst 
 \\[1mm] 
             & &\ \ +\ \xi(1 - \frac{1}{2} H_{13})H_{13}\otimes (\exmdst -\exmcst)
                 - 2 \xi E_{32} \otimes E_{23} e^{-3\widetilde{\sigma}} 
\\[1mm] 
             & & \ \ +\ 2 \xi E_{12} \otimes E_{21} e^{- \widetilde{\sigma}} 
 \\[1mm]
             & & \ \ - \ 4 \xi^2 (1 - \frac{1}{2} H_{13}) E_{12} \otimes 
                E_{23}e^{- \widetilde{\sigma}} (\exmdst - 2\exmcst) \\[1mm]               
             & & \ \ -\  4 \xi^2 E_{12} \otimes H_{13} E_{23}
e^{-3\widetilde{\sigma}} 
                 + 8 \xi^3 E_{12}^2 \otimes E_{23}^2 
e^{-6\widetilde{\sigma}} , 
                 \\[2mm]
\dele^{\rm can}(E_{21})&=& E_{21} \otimes  e^{- \widetilde{\sigma}} 
+ 1 \otimes E_{21} 
                 + \xi (H_{12} 
                 - H_{23}) \otimes E_{23} \exmdst. 
\end{array}
\ee 
Note that the deformation parameter $\xi$ and $\widetilde{\sigma} 
=\frac 12\ln (1+2\xi E)$ 
had been introduced here to make the correlations with the previous 
results more transparent.


\begin{thebibliography}{99}

\bibitem{D2}  V. G. Drinfeld, Leningrad Math. J. {\bf 1}, 1419 (1990).

\bibitem{D3}  V. G. Drinfeld, DAN USSR {\bf 273}, 531 (1983).

\bibitem{KUL1}  P. P. Kulish and A. A. Stolin, Czech. J. Phys. {\bf 47}, 123
(1997).

\bibitem{VLA} A. A. Vladimirov, Mod. Phys. Lett. {\bf A8}, 2573 (1993) 
(hep-th/9401101).

\bibitem{VAL} A. Ballesteros, F. J. Herranz, M. A. del Olmo, C. M.
Pere\~na and M. Santander, J. Phys. A: Math. Gen. {\bf 28}, 7113 (1995).

\bibitem{DRIN}  V. G. Drinfeld, ``Quantum groups", in {\em Proc. Int.
Congress of Mathematicians, Berkeley, 1986}, {\bf 1}.
Ed. A. V. Gleason (AMS, Providence, 1987).

\bibitem{OGIEV}  O. V. Ogievetsky,  Suppl. Rendiconti Cir. Math. Palermo,
Serie II {\bf  37}, 185 (1993) (preprint MPI-Ph/92-99, Munich (1992)).

\bibitem{GER} M.  Gerstenhaber, A. Giaquinto and S. D. Schak, Israel Mathem. 
Conference Proceedings, Vol. 7, 45 (1993).

\bibitem{KLM} P. P. Kulish, V. D. Lyakhovsky and A. I. Mudrov, ``Extended
jordanian twists for Lie algebras", math.QA/9806014 (submitted to 
J. Math. Phys).

\bibitem{RSTS}  N. Yu. Reshetikhin and M. A. Semenov-Tian-Shansky, J.
Geom. Phys. {\bf 5}, 533 (1988).

\bibitem{V-M} V. D. Lyakhovsky and M. A. del Olmo, ``Peripheric extended 
twists", (submitted to J. Phys. A) math.QA/9811153.  

\bibitem{JIMB}  M. Jimbo, Lett. Math. Phys. {\bf 10}, 63 (1985);
{\bf 11}, 247 (1986).

\bibitem{RES}  N. Yu. Reshetikhin, Lett. Math. Phys. {\bf 20},
331 (1990).

\bibitem{KM}  P. P. Kulish and A. I. Mudrov, ``Universal $R$--matrix for
esoteric quantum group", to be published in Lett. Math. Phys.
(math.QA/9804006).

\bibitem{KL} P. P. Kulish and V. D. Lyakhovsky, ``Classical and quantum
duality  in jordanian quantizations", Czech, J. Phys. {\bf 48}, 1415,
(1998),  (math.QA/9807122).

\bibitem{STO} A. Stolin, Math. Scand. {\bf 69}, 81 (1991).

\bibitem{ETI}  P. Etingof and D. Kazhdan, Selecta Math. {\bf 2}, 
1 (1996), ( q-alg/9510020).

\end{thebibliography}
\end{document}